\documentclass[a4paper]{amsart}
\usepackage{amssymb}
\usepackage{xypic}
\usepackage{graphicx}

\newtheorem{theorem}{Theorem}[section]

\newtheorem{cor}[theorem]{Corollary}

\theoremstyle{definition}

\theoremstyle{remark}

\numberwithin{equation}{section}

\newcommand{\C}{\mathbb{C}}

\DeclareMathOperator{\id}{Id}

\title{Equivariantly homeomorphic quasitoric manifolds are diffeomorphic}

\author{Michael Wiemeler}
\address{Institut f\"ur Mathematik\\ Universit\"at Augsburg\\D-86135 Augsburg\\Germany}
\email{michael.wiemeler@math.uni-augsburg.de}
\thanks{The research for this paper was partially supported by DFG grant HA 3160/6-1}



\begin{document}
\begin{abstract}
 In this note we prove that equivariantly homeomorphic quasitoric manifolds are diffeomorphic.
As a consequence we show that up to finite ambiguity the diffeomorphism type of certain quasitoric manifolds \(M\) is determined by their cohomology rings and first Pontrjagin classes.
\end{abstract}

\maketitle

%-----------------------------------------------------------------------
% End of amsart.template
%-----------------------------------------------------------------------

\section{Introduction}

A torus manifold is a \(2n\)-dimensional manifold with a smooth effective action of an \(n\)-dimensional torus \(T\) such that \(M^T\) is non-empty.
A torus manifold is called locally standard, if each point \(x\in M\) has an invariant neighborhood which is weakly equivariantly diffeomorphic to an open invariant subset of \(\C^n\), where \(T\) acts on \(\C^n\) by componentwise multiplication.
If this condition is satisfied the orbit space \(M/T\) is a smooth manifold with corners.

We call \(M\) quasitoric, if \(M/T\) is face-preserving homeomorphic to a simple convex polytope.
Examples of quasitoric manifolds are given by \(2n\)-dimensional symplectic manifolds with Hamiltonian actions of \(n\)-dimensional tori.

Quasitoric manifolds were introduced by Davis and Januszkiewicz \cite{MR1104531}. They showed that one can classify quasitoric manifolds up to equivariant homeomorphism by the combinatorial type of the orbit space and information on the isotropy groups of the torus action.

In \cite{MR3030690} we classified quasitoric manifolds up to equivariant diffeomorphism.
It turned out that there are equivariantly homeomorphic quasitoric manifolds which are not equivariantly diffeomorphic, if and only if there are exotic smooth structures on the four-dimensional disc.

In this note we address the question of classifying quasitoric manifolds up to non-equivariant diffeomorphism.
Our main result is as follows:

\begin{theorem}
\label{sec:equiv-home-quas}
  Let \(M\) and \(M'\) be two quasitoric manifolds which are equivariantly homeomorphic.
  Then \(M\) and \(M'\) are diffeomorphic.
\end{theorem}

Metaftsis and Prassidis \cite{metaftsis13:_topol} have shown that if a torus manifold \(M_1\) is equivariantly homotopy equivalent to a quasitoric manifold \(M_2\), then \(M_1\) and \(M_2\) are equivariantly homeomorphic. In particular, \(M_1\) is a quasitoric manifold.
Therefore we get the following corollary from Theorem~\ref{sec:equiv-home-quas}.

\begin{cor}
\label{sec:equiv-home-quas-1}
  Let \(M_1\) be a torus manifold and \(M_2\) be a quasitoric manifold.
  If \(M_1\) and \(M_2\) are equivariantly homotopy equivalent, then they are non-equivariantly diffeomorphic.
\end{cor}

In \cite{MR3030690} equivariantly homotopy equivalent torus manifolds were constructed which are not diffeomorphic.
The above corollary shows that such a construction can not be done in the category of quasitoric manifolds.

The question if a quasitoric manifold is determined by its cohomology ring is called the rigidity problem for quasitoric manifolds.
We can also prove the following theorem which gives some evidence for an answer to this question.

\begin{theorem}
\label{sec:introduction}
  Let \(M\) be a quasitoric manifold such that the intersection of any two facets of \(P=M/T\) is non-empty. Then the diffeomorphism type of \(M\) is determined up to finite ambiguity by its cohomology ring \(H^*(M;\mathbb{Z})\) and its first Pontrjagin class \(p_1(M)\in H^4(M;\mathbb{Z})\).
\end{theorem}

Note here that by a result of Panov and Ray \cite{MR2428364} all quasitoric manifolds are formal in the sense of rational homotopy theory.
Therefore by an argument of Sullivan \cite{MR0646078} (see also \cite{MR1113760}) the diffeomorphism type of a quasitoric manifold is determined up to finite ambiguity by its cohomology ring and its total Pontrjagin class.

A generalized Bott manifold is a manifold \(X_k\) such that there is a sequence of fibrations
\begin{equation*}
  X_k\rightarrow X_{k-1}\rightarrow \dots \rightarrow X_1\rightarrow X_0.
\end{equation*}
 where each \(X_i\), \(i>0\), is the total space of the projectivization of a sum of \(n_i+1\) line bundles over \(X_{i-1}\) and \(X_0\) is a point.
Note that the class of manifolds to which Theorem~\ref{sec:introduction} applies contains all the generalized Bott manifolds for which all \(n_i\) are greater than one.

A generalized Bott manifold is called Bott manifold if all \(n_i\) are equal to one.
Recently, it has been shown by Choi, Masuda and Murai \cite{zbMATH06442386}, that an isomorphism between the cohomology rings of two Bott manifolds preserves their Pontrjagin classes.
Therefore the diffeomorphism type of such a manifold is determined up to finite ambiguity by its cohomology ring.
For an overview of results concerning the rigidity of quasitoric manifolds and related questions see the survey \cite{MR2962979}.

Using \cite[Theorem 3.6]{wiemeler15:_torus} one can also prove a version of Theorem~\ref{sec:introduction} for simply connected torus manifolds whose cohomology ring is generated in degree two.
In this version the homeomorphism type of such a manifold is determined up to finite ambiguity by the cohomology ring and the first Pontrjagin class.

At the end of this introduction we give more concrete examples of applications of Theorem~\ref{sec:equiv-home-quas}.

In \cite[Corollary 2.7]{MR2885534} it has been shown that a quasitoric manifold \(M\) with the cohomology of a product \(\prod_i\C P^{n_i}\), \(n_i>1\), is weakly equivariantly homeomorphic to \(\prod_i\C P^{n_i}\) with the standard linear action.
Therefore such an \(M\) is also diffeomorphic to \(\prod_i \C P^{n_i}\).

Moreover, by the proof of \cite[Corollary 2.6]{MR2885534}, a quasitoric manifold with the same cohomology as \(\C P^n\#\C P^n\) is diffeomorphic to \(\C P^n\#\C P^n\).

I would like to thank the anonymous referee for several comments which helped to improve the presentation of this paper.

\section{Background on the classification of quasitoric manifolds}
\label{sec:background}

Let \(M\) be a quasitoric manifold, \(P=M/T\), \(\lambda\) the map which assigns to an orbit in \(M\) the isotropy group of this orbit.
The map \(\lambda\) is constant on the relative interiors of the faces of \(P\).
Therefore we may think of \(\lambda\) as defined on the set of faces of \(P\).
One easily sees that \(\lambda\) is uniquely determined by its values on the facets of \(P\).

It has been shown by Davis and Januszkiewicz \cite{MR1104531} that \(M\) is equivariantly homeomorphic to a standard model:
\begin{equation*}
  P\times T/\sim,
\end{equation*}
where \((x,t)\sim (x',t')\) if and only if \(x=x'\) and \(t^{-1}t'\in\lambda(x)\).
Therefore the equivariant homeomorphism type of \(M\) is determined by purely combinatorial data.

We can also construct a quasitoric manifold over the polytope \(P\) from a map
\[\lambda:\{\text{facets of } P\}\rightarrow \{S^1-\text{subgroups of } T\},\]
satisfying the following condition:
Whenever the facets \(F_1,\dots,F_n\) of \(P\) have non-trivial intersection, the natural map \(\lambda(F_1)\times\dots\times \lambda(F_n)\rightarrow T\) is an isomorphism.

Hence, we have a one-to-one correspondence between quasitoric manifolds (up to equivariant homeomorphism) and pairs \((P,\lambda)\), where \(P\) is a simple convex polytope and \(\lambda\) is a map satisfying the above conditions.

In \cite{MR3030690} we showed that the equivariant differentiable structures on a quasitoric manifold \(M\) (up to equivariant diffeomorphism) are in one-to-one correspondence to the smooth structures on \(P=M/T\) (up to diffeomorphism).
Therefore up to equivariant diffeomorphism quasitoric manifolds are in one-to-one correspondence to pairs \((X,\lambda)\) where \(X\) is a smooth manifold with corners which is face-preserving homeomorphic to a simple convex polytope and a map \(\lambda\) satisfying the above conditions.

Note that such an \(X\) is diffeomorphic to a simple convex polytope if and only if all its four-dimensional faces are diffeomorphic after smoothing the corners to standard discs \(D^4\).

The cohomology ring of a quasitoric manifold was described by Davis and Janusz\-kiewicz.
It is generated by the Poincar\'e-duals \(u_1,\dots,u_m\) of the characteristic submanifolds of \(M\).
A characteristic submanifold is the preimage of a facet of the orbit polytope under the orbit map.
They have codimension two.
The \(u_i\) are subject to two types of relations:
\begin{enumerate}
\item There are linear relations between the \(u_i\) depending on the map \(\lambda\);
\item A product \(u_{i_1}\dots u_{i_k}\) is zero if and only if the intersection of the characteristic submanifolds \(M_{i_1},\dots,M_{i_k}\) is empty.
\end{enumerate}

Therefore if the intersection of any two facets of the orbit polytope \(P\) of \(M\) is non-empty, then there are no relations of degree four between basis elements of \(H^2(M;\mathbb{Z})\).
Moreover, all \(u_i\) are non-zero.

The Pontrjagin classes of \(M\) can be computed as follows:
\begin{equation*}
  p(M)=\prod_{i=1}^m(1+u_i^2).
\end{equation*}

\section{The proof of Theorem \ref{sec:equiv-home-quas}}
\label{sec:proof}
  In this section we prove Theorem~\ref{sec:equiv-home-quas}.

  By Corollary~6.4 of \cite{MR3030690}, we may assume that \(M'\) and \(M\) have dimension at least eight.
By Corollary 5.7 of \cite{MR3030690}, two equivariantly homeomorphic quasitoric manifolds whose orbit spaces are diffeomorphic to simple convex polytopes are equivariantly diffeomorphic.
  Therefore it is sufficient to show that we can change the torus actions on \(M'\) and \(M\) (without changing their equivariant homeomorphism type) in such a way that \(M'/T\) and \(M/T\) become diffeomorphic to the simple convex polytope \(P\) which is face-preserving homeomorphic to \(M'/T\) and \(M/T\).
  By the symmetry of the problem, it suffices to do this construction for \(M'\).
  As noted in Section~\ref{sec:background} the orbit space \(M'/T\) is diffeomorphic to \(P\) if all four-dimensional faces of \(M'/T\) are diffeomorphic after smoothing the corners to standard discs.

  First assume that \(\dim M =2n \geq 10\).
  In this case the proof combines ideas from the proofs of Theorems~5.5 and 7.1 of \cite{MR3030690} (see also the discussion in Section~3 of \cite{wiemeler15:_torus}).
  Let \(X\) be a four dimensional face of \(M'/T\) which is not diffeomorphic after smoothing the corners to a standard disc.
  Then \(\partial X\) is diffeomorphic to \(S^3\).
  Hence, \(X\cup_{S^3}D^4\) is a four-dimensional homotopy sphere and therefore bounds a contractible five dimensional manifold \(S\) \cite[Theorem 3]{MR0253347}.

  Now form \[\alpha(M'/T,X)=(M'/T-X\times \Delta^{n-4})\cup_{X\times \Delta^{n-5}} (S\times \Delta^{n-5}\cup_{D^4\times \Delta^{n-5}}D^4\times \Delta^{n-4}),\]
where \(\Delta^m\) denotes the \(m\)-dimensional simplex.
Loosely speaking one can think of \(\alpha(M'/T,X)\) as constructed from \(M'/T\) by first removing the face \(X\) and then replacing it by a face which is diffeomorphic to \(D^4\) after smoothing the corners.

  Then \(\alpha(M'/T,X)\) is a nice manifold with corners
  such that all its faces are contractible.

  To see this last statement first note that there is a one-to-one correspondence \(\phi\) between the faces of \(M'/T\) and \(\alpha(M'/T,X)\).
  This one-to-one correspondence is defined as follows: If \(Y\) is a face of \(M'/T\) which is not contained in \(X\), then the face \(\phi(Y)\) of \(\alpha(M'/T,X)\) is the unique face of \(\alpha(M'/T,X)\), such that 
\[Y\cap (M'/T-X\times \Delta^{n-4})=\phi(Y)\cap (M'/T-X\times \Delta^{n-4}).\]
This correspondence can be extended to the set of all faces of \(M'/T\) in such a way that it is compatible with intersections of faces, i.e. such that if \(Y=\bigcap_i Y_i\), then we have \(\phi(Y)=\bigcap_i \phi(Y_i)\) for faces \(Y\) and \(Y_i\) of \(M'/T\).

  Now the following holds:
  If \(Y\) is a face of \(M'/T\) which does not contain \(X\), then \(\phi(Y)\) is diffeomorphic to \(Y\).
  Moreover, if \(Y\) contains \(X\) properly, then the face \(\phi(Y)\) can be constructed by gluing \(Y\) and a contractible manifold along a contractible subset of the boundary of \(Y\).
  So \(\phi(Y)\) is contractible.
  The face which corresponds to \(X\) is diffeomorphic after smoothing the corners to \(D^4\).

  Hence, by Lemma 3.2 of \cite{wiemeler15:_torus}, \(\alpha(M'/T,X)\) 
  is face-preserving homeomorphic to \(M'/T\).
  
  Let \(M''\) be the quasitoric manifold corresponding to \((\alpha(M'/T,X),\lambda)\). We show that \(M''\) and \(M'\) are diffeomorphic.

  Let \(\Sigma_S= \partial (S\times D^{2(n-4)})\times T^4\).
  Then, by the h-cobordism Theorem applied to \((S\times D^{2(n-4)})- \dot{D}^{2n-4}\), where \(\dot{D}^{2n-4}\) is a tubular neighborhood of an interior point of \(S\times D^{2(n-4)}\), \(\Sigma_S\) is diffeomorphic to \(S^{2n-4}\times T^{4}\).
By a similar argument, \(X\times D^{2(n-4)}\) is diffeomorphic to \(D^{2n-4}\).
  By Schoenfies' Theorem we can assume that \(X\times D^{2(n-4)}\) is embedded in \(S^{2n-4}\) as the upper hemisphere. 
  Moreover, \[M''=(M'-(\dot{X}\times D^{2(n-4)}\times T^4))\cup_{\partial (\dot{X}\times D^{2(n-4)})\times T^4}(\Sigma_S-(\dot{X}\times D^{2(n-4)})\times T^4),\]
where \(\dot{X}\) denotes \(X\) with a small collar of its boundary removed.
  Therefore \(M''\) and \(M'\) are diffeomorphic.
  The claim now follows by induction on the number of four-dimensional faces of \(M'/T\) which are not diffeomorphic to \(D^4\) after smoothing the corners.
  
  Now assume that \(\dim M=8\).
  Let \(D_1\) be the orbit space of \(M\) with a small collar of its boundary removed.
  Similarly, define \(D_2\) to be the orbit space of \(M'\) with a small collar of its boundary removed.
  Then \(D_1\) and \(D_2\) are homeomorphic to \(D^4\).
  Moreover, by the proof of Theorem 5.4 of \cite{MR3030690}, \(M/T-D_1\) and \(M'/T-D_2\) are diffeomorphic.
  Therefore, by the proof of Theorem 5.6 of \cite{MR3030690}, there is an eight-dimensional \(T\)-manifold \(N\)  with boundary such that
  \begin{align*}
    M&= D_1\times T^4\cup_{g_1} N& M'&=D_2\times T^4\cup_{g_2}N,
  \end{align*}
where \(g_i:S^3\times T^4=\partial D_i\times T^4 \rightarrow \partial N\) are equivariant diffeomorphisms.

Note that the orbit spaces of \(M\) and \(M'\) are diffeomorphic after smoothing the corners to \(D_1\) and \(D_2\), respectively.
Therefore, by the remarks at the beginning of this section, it suffices to consider the case that \(D_1\) is diffeomorphic to \(D^4\).
Moreover, since \(N/T\) is diffeomorphic after smoothing the corners to \(S^3\times I\), there is a homeomorphism \(f: N/T\rightarrow S^3\times I\)
which is a diffeomorphism outside a small neighborhood of the codimension-two faces of \(N/T\) such that \(f((\partial N)/T)=S^3\times \{0\}\).
Moreover, because every diffeomorphism of \(S^3\times\{0\}\) extends to a diffeomorphism of \(S^3\times I\), we may assume that \(f^{-1}|_{S^3\times \{0\}}\) is the map which is induced by \(g_1\) on the orbit space.

Furthermore, there is an \(T^4\)-equivariant homeomorphism
\begin{equation*}
  f':N\rightarrow (N/T^4\times T^4)/\sim,
\end{equation*}
where \((x_1,t_1)\sim(x_2,t_2)\) if and only if \(x_1=x_2\) and \(t_1t_2^{-1}\in \lambda(x_1)\).
Here \(\lambda(x_1)\) denotes the isotropy group of the orbit \(x_1\).

Because \(g_2^{-1}\circ g_1\) is \(T^4\)-equivariant, there is a diffeomorphism \(g':S^3\rightarrow S^3\) and a map \(g'':S^3\rightarrow T^4\) such that, for \((x,y)\in S^3\times T^4\),
\begin{equation*}
  g_2^{-1}\circ g_1(x,y)=(g'(x),g''(x)y).
\end{equation*}
Since \(\pi_3(T^4)=0\), we may assume that \(g''(x)=1\) for all \(x\in S^3\).
Hence \(g_2^{-1}\circ g_1\) is given by \(g'\times \id_{T^4}\).
We write \(T^4=T^2\times T'^2\) such that \(\lambda(F)\subset T'^2\) for a facet \(F\) of \(N/T\).

The Whitehead torsion group of a free abelian group is trivial \cite{MR0174605}.
Therefore it follows from the s-cobordism theorem applied to \((D_i-\dot{D}^4)\times T^2\), where \(\dot{D}^4\) is a tubular neighborhood of an interior point of \(D_i\), that \(D_i\times T^2\), \(i=1,2\), is diffeomorphic to \(D^4\times T^2\).
Let \(h:D_2\times T^2\rightarrow D_1\times T^2\) be a diffeomorphism.

~

\emph{Claim:} \(h\circ (g'\times \id_{T^2}):S^3\times T^2\rightarrow S^3\times T^2\) is isotopic to a diffeomorphism of the form \(h_1\circ h_2\), where \(h_1:S^3\times T^2\rightarrow S^3\times T^2\) is a diffeomorphism of the form
\begin{equation}
\label{eq:2}
  (x,t)\mapsto (\tilde{h}_1(t)x,t)
\end{equation}
for some map \(\tilde{h}_1:T^2\rightarrow O(3)\subset O(4)\)
and \(h_2:S^3\times T^2\rightarrow S^3\times T^2\) is a diffeomorphism such that \(h_2|_{D^3\times T^2}=\id_{D^3\times T^2}\) for some disc \(D^3\subset S^3\).

~

It follows from Haefliger's Theorem \cite{MR0145538} that \(h\circ (g'\times \id_{T^2})\) is isotopic to a diffeomorphism \(h'\) which is the identity on \(\{pt\}\times T^2\) for the north pole \(pt\in S^3\).
Now choose a tubular neighborhood of \(\{pt\}\times T^2\) of the form \(D^3\times T^2\), where \(D^3\subset S^3\) is the upper hemisphere.
Then \(h'(D^3\times T^2)\) is also a tubular neighborhood of \(\{pt\}\times T^2\).
Therefore \(h'\) can be isotoped to a diffeomorphism \(h''\) such that \(h''|_{D^3\times T^2}\) is given by
\begin{equation*}
  (x,t)\mapsto (\tilde{h}_1(t)x,t)
\end{equation*}
with \(\tilde{h}_1\) as above.

Now \(h_1\) is given by (\ref{eq:2}) and \(h_2=h_1^{-1}\circ h''\).
This proves the claim.

Since \(D_1\) is diffeomorphic to \(D^4\), \(h_1\) can be extended to a diffeomorphism of \(D_1\times T^2\).
Therefore we have
\begin{align*}
  M'&=D_2\times T^4 \cup_{g_2} N\\
  &=D_1\times T^4\cup_{g_2\circ (h^{-1}\times \id_{T'^2})} N\\
  &=D_1\times T^4\cup_{g_1\circ (g'^{-1}\times \id_{T^4})\circ (h^{-1}\times \id_{T'^2})} N\\
  &=D_1\times T^4\cup_{g_1\circ (h_2^{-1}\times \id_{T'^2})} N.
\end{align*}
Since \(M=D_1\times T^4\cup_{g_1} N\), it now follows that the identity map on \(D_1 \times T^4\) extends to a diffeomorphism \(M \rightarrow M'\) if
\(g_1 \circ (h^{-1}_2\times Id_{T'^2}) \circ g_1^{-1} : \partial N \rightarrow \partial N\) extends to a diffeomorphism from \(N\) to itself.

\begin{figure}
  \centering
  \includegraphics{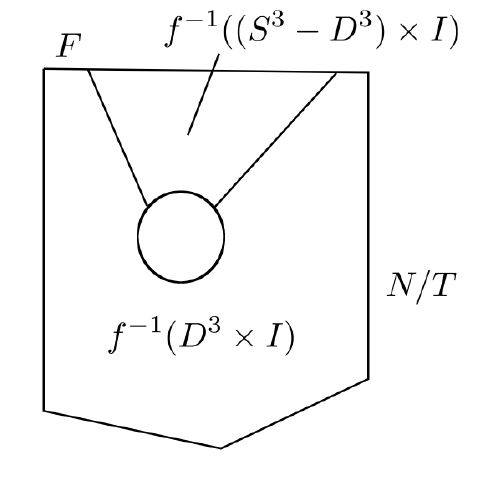}
  \caption{The orbit space of $N$ and the homeomorphism $f$.}
  \label{fig:1}
\end{figure}

Now consider \(f^{-1}((S^3-D^3)\times I)\).
It intersects the boundary of \(N/T\) in \(f^{-1}((S^3-D^3)\times\{0\})\) and in \(f^{-1}((S^3-D^3)\times \{1\})\).

Since \(S^3-D^3\) can be identified with a tubular neighborhood of a point in \(S^3\), we may assume that \(S^3-D^3\) is so small that \(f^{-1}((S^3-D^3)\times \{1\})\) is contained in the interior of the facet \(F\) of \(N/T\).
(See Figure \ref{fig:1}.)
Moreover, since \(f\) is a diffeomorphism outside a small neighborhood of the codimension-two faces of \(N/T\), we may assume that the restriction of \(f^{-1}\) to \((S^3-D^3)\times I\) is a diffeomorphism onto its image.

Hence, by \cite[Theorem 6.3, p. 331]{0246.57017}, there is an \(T^4\)-equivariant diffeomorphism
\begin{equation*}
  f'':\pi^{-1}(f^{-1}((S^3-D^3)\times I))\rightarrow (S^3-D^3)\times T^2\times (T'^2\times_{\lambda(F)}D^2)
\end{equation*}
which induces \(f|_{f^{-1}((S^3-D^3)\times I)}\) on the orbit space.

Here \(\pi\) denotes the orbit map.
Since \(g_1\circ (h_2^{-1}\times \id_{T'^2})\circ g_1^{-1}\) is \(T'^2\)-equivariant
and \(\lambda(F)\subset T'^2\),
the diffeomorphism
\begin{align*}
  ((S^3-D^3)\times T^2\times T'^2)\times D^2&\rightarrow ((S^3-D^3)\times T^2\times T'^2)\times D^2\\
  (x,y)&\mapsto ((f''\circ g_1\circ (h_2^{-1}\times \id_{T'^2})\circ g_1^{-1}\circ f''^{-1})(x), y)
\end{align*}
induces a diffeomorphism of \(\pi^{-1}(f^{-1}((S^3-D^3)\times I))\) which extends \(g_1\circ (h_2^{-1}\times \id_{T'^2})\circ g_1^{-1}|_{\pi^{-1}(f^{-1}((S^3-D^3)\times\{0\}))}\).
 
Moreover, since \(g_1\circ (h_2^{-1}\times \id_{T'^2})\circ g_1^{-1}|_{\pi^{-1}(f^{-1}((D^3)\times \{0\}))}\) is the identity, it  can be extended by the identity on \(\pi^{-1}(f^{-1}(D^3\times I))\).
These two extensions of \(g_1\circ (h_2^{-1}\times \id_{T'^2})\circ g_1^{-1}\) fit together to form a diffeomorphism of \(N\).

From this it follows that the identity on \(D_1\times T^4\) can be extended to a diffeomorphism \(M\rightarrow M'\).
This proves the theorem.

\section{The proof of Theorem~\ref{sec:introduction}}

In this section we prove Theorem~\ref{sec:introduction}.
  By Theorem~\ref{sec:equiv-home-quas} and Theorem 2.2 of \cite{MR2885534}, it is sufficient to show that only finitely many classes in \(H^2(M;\mathbb{Z})\) can be realized by the Poincar\'e-duals of the characteristic submanifolds of a quasitoric manifold \(M'\) with \(H^*(M;\mathbb{Z})=H^*(M';\mathbb{Z})\) and \(p_1(M)=p_1(M')\).
Let \(v_1,\dots,v_k\) be a basis of \(H^2(M;\mathbb{Z})\).
By the assumption on the combinatorial type of \(P\), a basis of \(H^4(M;\mathbb{Z})\) is then given by \(v_iv_j\), \(i,j=1,\dots,k\).

Let \(u_1,\dots,u_m\) the Poincar\'e-duals of the characteristic submanifolds of \(M'\).
Note that by the description of the cohomology ring of \(M'\) given in Section~\ref{sec:background} we have \(m=b_2(M')+\frac{1}{2}\dim M'=k+n\).
 Let, moreover, \(\mu_{ij}\in \mathbb{Z}\) such that \(u_i=\sum_{j=1}^k\mu_{ij} v_j\).

Then we have
\begin{equation*}
  p_1(M)=p_1(M')=\sum_{i=1}^m u_i^2=\sum_{j=1}^k \left(\sum_{i=1}^m \mu_{ij}^2\right)v_j^2 + \sum_{i\neq j} \nu_{ij}v_iv_j
\end{equation*}
with some \(\nu_{ij}\in \mathbb{Z}\).
Therefore there is a constant \(C\) which depends on \(p_1(M)\) and the \(v_j\) such that \(|\mu_{ij}|<C\) for all \(i,j\).
This means that the \(u_i\) are contained in a finite set.
This proves the theorem.

\bibliography{equi_homeo_qt}{}
\bibliographystyle{amsplain}
\end{document}